\documentclass[10pt,oneside]{article}
\input amssymb.sty
\usepackage[cp1250]{inputenc}
\usepackage[T1]{fontenc}

\begin{document}
\newcommand{\qed}{\rule{1.5mm}{1.5mm}}
\newcommand{\proof}{\textit{Proof. }}
\newcommand{\ccon}{\rightarrowtail}
\newtheorem{theorem}{Theorem}[section]
\newtheorem{lemma}[theorem]{Lemma}
\newtheorem{remark}[theorem]{Remark}
\newtheorem{example}[theorem]{Example}
\newtheorem{corollary}[theorem]{Corollary}
\newtheorem{proposition}[theorem]{Proposition}
\newtheorem{claim}[theorem]{Claim}
\begin{center}
\vspace*{0mm}
{\LARGE\textbf{Approximation of sets defined by polynomials with holomorphic coefficients}\vspace*{3mm}}\\
{\large Marcin Bilski\footnote{M. Bilski: Institute of
Mathematics, Jagiellonian University, Reymonta 4, 30-059 Krak\'ow,
Poland. e-mail: Marcin.Bilski@im.uj.edu.pl\vspace*{1mm}\\
Research partially supported by the Polish Ministry of Science and
Higher Education.}}\vspace*{10mm}\\
\end{center}
\begin{abstract}Let $X$ be an analytic set defined by polynomials whose
coefficients $a_1,\ldots,a_s$ are holomorphic functions. We
formulate conditions such that for all sequences
$\{a_{1,\nu}\},\ldots,\{a_{s,\nu}\}$ of holomorphic functions
converging locally uniformly to $a_1,\ldots,a_s$ respectively the
following holds true. If $a_{1,\nu},\ldots,a_{s,\nu}$ satisfy the
conditions then the sequence of the sets $\{X_{\nu}\}$ obtained by
replacing $a_j$'s by $a_{j,\nu}$'s in the polynomials, converge to
$X.$\vspace*{2mm}\\
\textbf{Keywords } Analytic set, Nash set,
approximation\vspace*{0mm}\\
\textbf{MSC (2000):} 32C25
\end{abstract}
%+++++++++++++++++++++++++++++++++++++++++++++++++++++++++++++++++
\section{Introduction and main results}
\label{intro} The problem of approximating analytic objects by
simpler algebraic ones with similar properties appears in many
contexts of complex geometry and has\linebreak attracted the
attention of several mathematicians (see \cite{Ar2}, \cite{Ar3},
\cite{BoK}, \cite{BMT1}, \cite{DLS}, \cite{Fo}, \cite{Ku},
\cite{Lem}, \cite{TT1}, \cite{TT2}, \cite{TT3}). In the present
paper we concern the problem in the case where the approximated
objects are complex analytic sets whereas the approximating ones
are complex Nash sets (see Section~ \ref{prelnash}). The
approximation is expressed by means of the convergence of
holomorphic chains (for the definition see Section~
\ref{holchai}).

For sets with proper projection the existence of such
approximation was discussed in \cite{B1}, \cite{B2}. In a
subsequent paper \cite{B3} it was proved that the order of
tangency of a limit set and the approximating sets can be
arbitrarily high. The first results on approximation of analytic
sets by higher order tangent algebraic varieties are due to
R.~W.~Braun, R. Meise and B. A. Taylor \cite{BMT1} with
applications in \cite{BMT2}.

Both in \cite{B2} and in \cite{B3} analytic sets are represented
as mappings defined on an open subset of $\mathbf{C}^n$ with
values in an appropriate symmetric power of $\mathbf{C}^m.$
However, in many cases such sets are defined by systems of
equations which in general carry more information than the sets
themselves. Therefore it is natural to look for approximations of
the functions appearing in the equations. Throughout this paper we
restrict our attention to the case where the description is given
by a system of polynomials with holomorphic coefficients whereas
the approximated set is with proper projection onto an appropriate
affine space. Our aim is to show how to approximate the
coefficients of the polynomials to obtain Nash approximations of
the set.

If the number of the functions describing the analytic set $X$ is
equal to the codimension of $X$ then it is sufficient to take
generic approximations of the coefficients in order to get local
uniform approximation of $X.$ Such approach clearly does not work
in the case of a non-complete intersection as it leads to sets of
dimensions strictly smaller than the dimension of $X.$ Yet, it is
natural to expect that there are algebraic relations satisfied by
the coefficients such that if the approximating coefficients also
satisfy the relations then the original polynomials with these new
coefficients define appropriate approximations.

Before stating the main result let us recall that for any analytic
set $Y$ by $Y_{(n)}$ we denote the union of all $n$-dimensional
irreducible components of $Y.$

Let $U\subset\mathbf{C}^n$ be a domain. Abbreviate
$v=(v_1,\ldots,v_p),$ $z=(z_1,\ldots,z_m).$ Assuming the notation
of Section~\ref{prel} and treating analytic sets as holomorphic
chains with components of multiplicity one we prove
\begin{theorem}\label{main}Let $q_1,\ldots,q_s\in\mathbf{C}[v,z],$
for some $s\in\mathbf{N},$ and let $H:U\rightarrow\mathbf{C}^p$ be
a holomorphic mapping. Assume that
$$X=\{(x,z)\in U\times\mathbf{C}^m: q_i(H(x),z)=0,
i=1,\ldots,s\}$$ is an analytic set of pure dimension $n$ with
proper projection onto $U.$ Then there is an algebraic subvariety
$F$ of $\mathbf{C}^p$ with $H(U)\subset F$ such that for every
sequence $\{H_{\nu}:U\rightarrow F\}$ of holomorphic mappings
converging locally uniformly to $H$ the following holds. The
sequence $\{X_{\nu}\},$ where
$$X_{\nu}=\{(x,z)\in U\times\mathbf{C}^m:
q_i(H_{\nu}(x),z)=0, i=1,\ldots,s\},$$ converges to $X$ locally
uniformly and the sequence $\{(X_{\nu})_{(n)}\}$ converges to $X$
in the sense of holomorphic chains.
\end{theorem}

The following example shows that the sets from $\{X_{\nu}\}$ are
in general not purely dimensional:
\begin{example}
Define $X=\{(x,z)\in\mathbf{C}^2:zxe^x=0, z^2-zx=0\}.$ Then
$X=\{(x,z)\in\mathbf{C}^2:z=0\},$ therefore it is purely
$1$-dimensional. On the other hand, $\mathbf{C}^2\times\{1\}$ is
the smallest algebraic set in $\mathbf{C}^3$ containing the image
of the mapping $x\mapsto (-x,xe^x,1).$ By approximating this
mapping by $x\mapsto(-x,(x-\frac{1}{\nu})e^x,1)$ one obtains
$X_{\nu}=\{(x,z)\in\mathbf{C}^2: z(x-\frac{1}{\nu})e^x=0,
z^2-zx=0\}$ containing an isolated point
$(\frac{1}{\nu},\frac{1}{\nu}).$
\end{example}

Let $U$ be a connected Runge domain in $\mathbf{C}^n,$ let $X$ be
a purely $n$-dimen-\linebreak sional analytic subset of
$U\times\mathbf{C}^m$ with proper projection onto $U$ and let
$Q_1,\ldots,Q_{s}$ $\in\mathcal{O}(U)[z],$ for some
$s\in\mathbf{N},$ satisfy
$$X=\{(x,z)\in U\times\mathbf{C}^m:
Q_1(x,z)=\ldots=Q_{s}(x,z)=0\}.$$ (An example of such
$Q_1,\ldots,Q_s$ are the canonical defining functions for $X$ see
\cite{Wh}, \cite{Ch}).)

We check that combining Theorem \ref{main} with one of results of
L. Lempert (Theorem~3.2 from \cite{Lem}, see Theorem~\ref{lempth}
below) one obtains Nash approximations of $X$ by approximating its
holomorphic description by a Nash description. (Let us mention
that the proof of Theorem~\ref{lempth} is based on the affirmative
solution to the Artin's conjecture first presented in \cite{Po1},
\cite{Po2}, see also \cite{An}, \cite{Og}, \cite{Sp}.)

Let $H=(H_1,\ldots,H_{s})$ denote the holomorphic mapping defined
on $U$ where, for every $j\in\{1,\ldots,s\},$ $H_j$ is the mapping
whose components are all the non-zero coefficients of the
polynomial $Q_j;$ by $n_j$ denote the number of these
coefficients. More precisely, the components of $H_j$ are indexed
by $m$-tuples from some finite set $S_j\subset\mathbf{N}^m$ in
such a way that the component indexed by a fixed
$(\alpha_1,\ldots,\alpha_m)$ is the coefficient standing at the
monomial $z_1^{\alpha_1}\cdot\ldots\cdot z_m^{\alpha_m}$ in $Q_j.$

Let $F$ be the intersection of all algebraic subvarieties of
$\mathbf{C}^{(\sum_jn_j)}$ containing $H(U)$ and let $\tilde{U}$
be any open relatively compact subset of $U.$ Then $\tilde{U}$ is
contained in a polynomially convex compact subset of $U$ hence by
Theorem \ref{lempth} there exists a sequence
$\{H_{\nu}:\tilde{U}\rightarrow F\}$ of Nash mappings,
$H_{\nu}=(H_{1,\nu},\ldots, H_{s,\nu}),$ such that $\{H_{j,\nu}\}$
converges uniformly to $H_{j}|_{\tilde{U}},$ for every
$j=1,\ldots,s.$ Now let
$$X_{\nu}=\{(x,z)\in\tilde{U}\times\mathbf{C}^m:Q_{1,\nu}(x,z)=\ldots=Q_{s,\nu}(x,z)=0\},$$ where
$Q_{j,\nu}\in\mathcal{O}(\tilde{U})[z],$ for $j=1,\ldots,s,$ is
defined as follows. The coefficient of $Q_{j,\nu}$ standing at the
monomial $z_1^{\alpha_1}\cdot\ldots\cdot z_m^{\alpha_m}$ is the
component of $H_{j,\nu}$ indexed by $(\alpha_1,\ldots,\alpha_m)$
(if $(\alpha_1,\ldots,\alpha_m)\notin S_j$ then the coefficient
equals zero).

Finally, let $q_1,\ldots,q_s$ be the polynomials obtained from
$Q_1,\ldots,Q_s$ by replacing the holomorphic coefficients of the
latter polynomials by independent new variables. It is easy to see
that $q_1,\ldots,q_s$ together with the mapping $H$ satisfy the
hypotheses of Theorem~\ref{main}. Hence the sequence of Nash sets
$\{(X_{\nu})_{(n)}\},$ where $X_{\nu}$ defined in the previous
paragraph, converges to $X\cap(\tilde{U}\times\mathbf{C}^m)$ in
the sense of holomorphic chains. Thus we recover the main result
of \cite{B2}:
\begin{corollary}\label{conv} Let $X$ be a purely $n$-dimensional analytic subset of
$U\times\mathbf{C}^m$ with proper projection onto $U.$ Then for
every open set $\tilde{U}\subset\subset U$ there is a sequence
$\{X_{\nu}\}$ of purely $n$-dimensional Nash subsets of
$\tilde{U}\times\mathbf{C}^m$ converging to
$X\cap(\tilde{U}\times\mathbf{C}^m)$ in the sense of chains.
\end{corollary}

Every purely $n$-dimensional analytic set is locally with proper
projection onto an open subset of an $n$-dimensional affine space.
Hence, by Corollary \ref{conv} every analytic set can be locally
approximated by Nash ones. This fact is derived in a purely
geometric way (without using N\'eron desingularization standing
behind Theorem 3.2) in \cite{B4}.

Note that the convergence of positive chains appearing in this
paper is equivalent to the convergence of currents of integration
over the considered sets (see \cite{Lel}, \cite{Ch}). The
organization of this paper is as follows. In Section~ \ref{prel}
preliminary material is presented whereas Section~ \ref{mainproof}
contains the proof of Theorem \ref{main}.
%******************************************************************
\section{Preliminaries}
\label{prel}
%*******************************************************************************************
\subsection{Nash sets}\label{prelnash}
Let $\Omega$ be an open subset of $\mathbf{C}^n$ and let $f$ be a
holomorphic function on $\Omega.$ We say that $f$ is a Nash
function at $x_0\in\Omega$ if there exist an open neighborhood $U$
of $x_0$ and a polynomial
$P:\mathbf{C}^n\times\mathbf{C}\rightarrow\mathbf{C},$ $P\neq 0,$
such that $P(x,f(x))=0$ for $x\in U.$ A holomorphic function
defined on $\Omega$ is said to be a Nash function if it is a Nash
function at every point of $\Omega.$ A holomorphic mapping defined
on $\Omega$ with values in $\mathbf{C}^N$ is said to be a Nash
mapping if each of its components is a Nash function.

A subset $Y$ of an open set $\Omega\subset\mathbf{C}^n$ is said to
be a Nash subset of $\Omega$ if and only if for every
$y_0\in\Omega$ there exists a neighborhood $U$ of $y_0$ in
$\Omega$ and there exist Nash functions $f_1,\ldots,f_s$ on $U$
such that $$Y\cap U=\{x\in U: f_1(x)=\ldots=f_s(x)=0\}.$$

The fact from \cite{Tw} stated below explains the relation between
Nash and algebraic sets.
\begin{theorem}
Let $X$ be an irreducible Nash subset of an open set
$\Omega\subset\mathbf{C}^n.$ Then there exists an algebraic subset
$Y$ of $\mathbf{C}^n$ such that $X$ is an analytic irreducible
component of $Y\cap\Omega.$ Conversely, every analytic irreducible
component of $Y\cap\Omega$ is an irreducible Nash subset of
$\Omega.$
\end{theorem}
%**************************************************************************************************
\subsection{Convergence of closed sets and holomorphic chains}\label{holchai}
Let $U$ be an open subset in $\mathbf{C}^m.$ By a holomorphic
chain in $U$ we mean the formal sum $A=\sum_{j\in J}\alpha_jC_j,$
where $\alpha_j\neq 0$ for $j\in J$ are integers and
$\{C_j\}_{j\in J}$ is a locally finite family of pairwise distinct
irreducible analytic subsets of $U$ (see \cite{Tw2}, cp. also
\cite{Ba}, \cite{Ch}). The set $\bigcup_{j\in J}C_j$ is called the
support of $A$ and is denoted by $|A|$ whereas the sets $C_j$ are
called the components of $A$ with multiplicities $\alpha_j.$ The
chain $A$ is called positive if $\alpha_j>0$ for all $j\in J.$ If
all the components of $A$ have the same dimension $n$ then $A$
will be called an $n-$chain.

Below we introduce the convergence of holomorphic chains in $U$.
To do this we first need the notion of the local uniform
convergence of closed sets. Let $Y,Y_{\nu}$ be closed subsets of
$U$ for $\nu\in\mathbf{N}.$ We say that $\{Y_{\nu}\}$ converges to
$Y$ locally uniformly if:\vspace*{2mm}\\
(1l) for every $a\in Y$ there exists a sequence $\{a_{\nu}\}$ such
that $a_{\nu}\in Y_{\nu}$
and\linebreak\hspace*{7mm}$a_{\nu}\rightarrow a$
in the standard topology of $\mathbf{C}^m,$\\
\\
(2l)  for every compact subset $K$ of $U$ such that $K\cap
Y=\emptyset$ it holds $K\cap Y_{\nu}=\emptyset$\linebreak\hspace*{6.3mm}for almost all $\nu.$\vspace*{1mm}\\
Then we write $Y_{\nu}\rightarrow Y.$ For details concerning the
topology of local uniform convergence see \cite{Tw2}.

We say that a sequence $\{Z_{\nu}\}$ of positive $n$-chains
converges to a positive $n$-chain $Z$ if:\vspace*{2mm}\\
(1c) $|Z_{\nu}|\rightarrow |Z|,$\\
(2c) for each regular point $a$ of $|Z|$ and each submanifold $T$
of $U$ of dimension\linebreak\hspace*{7mm}$m-n$ transversal to
$|Z|$ at $a$ such that $\overline{T}$ is compact and
$|Z|\cap\overline{T}=\{a\},$\linebreak\hspace*{7mm}we have
$deg(Z_{\nu}\cdot T)=deg(Z\cdot T)$ for almost
all $\nu.$\vspace*{2mm}\\
Then we write $Z_{\nu}\ccon Z.$ (By $Z\cdot T$ we denote the
intersection product of $Z$ and $T$ (cf. \cite{Tw2}). Observe that
the chains $Z_{\nu}\cdot T$ and $Z\cdot T$ for sufficiently large
$\nu$ have finite supports and the degrees are well defined.
Recall that for a chain $A=\sum_{j=1}^d\alpha_j\{a_j\},$
$deg(A)=\sum_{j=1}^d\alpha_j$).

The following lemma from \cite{Tw2} will be useful to us.
\begin{lemma}\label{eqconv} Let $n\in\mathbf{N}$ and $Z,Z_{\nu},$
for $\nu\in\mathbf{N},$ be positive $n$-chains. If
$|Z_{\nu}|\rightarrow |Z|$ then the following conditions are
equivalent:\\
(1) $Z_{\nu}\ccon Z,$\\
(2) for each point $a$ from a given dense subset of $Reg(|Z|)$
there exists a sub-\linebreak\hspace*{5.5mm}manifold $T$ of $U$ of
dimension $m-n$ transversal to $|Z|$ at $a$ such that
$\overline{T}$\linebreak\hspace*{5.5mm}is compact,
$|Z|\cap\overline{T}=\{a\}$ and $deg(Z_{\nu}\cdot T)=deg(Z\cdot
T)$ for almost all $\nu.$
\end{lemma}

\subsection{Approximation of holomorphic mappings} In the proofs of Corollary \ref{conv}
we use the following theorem which is due to L.~Lempert (see
\cite{Lem}, Theorem 3.2).
\begin{theorem}
\label{lempth} Let $K$ be a holomorphically convex compact subset
of $\mathbf{C}^n$ and $f:K\rightarrow\mathbf{C}^k$ a holomorphic
mapping that satisfies a system of equations $Q(z,f(z))=0$ for
$z\in K.$ Here $Q$ is a Nash mapping from a neighborhood
$U\subset\mathbf{C}^n\times\mathbf{C}^k$ of the graph of $f$ into
some $\mathbf{C}^q.$ Then $f$ can be uniformly approximated by a
Nash mapping $F:K\rightarrow\mathbf{C}^k$ satisfying
$Q(z,F(z))=0.$
\end{theorem}
%***********************************************************************************************
\section{Proof of Theorem \ref{main}}\label{mainproof}
Denote $B_m(r)=\{z\in\mathbf{C}^m:||z||_{\mathbf{C}^m}<r\}$ and
recall $v=(v_1,\ldots,v_p).$ Let $U$ be a domain in
$\mathbf{C}^n.$ We prove the following
\begin{proposition}\label{ndimensional}Let $q_1,\ldots,q_s\in\mathbf{C}[v,z],$
for some $s\in\mathbf{N},$ and let $H:U\rightarrow\mathbf{C}^p$ be
a holomorphic mapping. Assume that
$$X=\{(x,z)\in U\times\mathbf{C}^m: q_i(H(x),z)=0,
i=1,\ldots,s\}$$ is an analytic set of pure dimension $n$ with
proper projection onto $U.$ Then there is an algebraic subvariety
$F$ of $\mathbf{C}^p$ with $H(U)\subset F$ such that for every
domain $\tilde{U}\subset\subset U$ and every sequence
$\{H_{\nu}:\tilde{U}\rightarrow F\}$ of holomorphic mappings
converging uniformly to $H$ on $\tilde{U}$ the following holds.
There is $r_0>0$ such that for every $r>r_0$ the sequence
$\{X_{\nu}\},$ where
$$X_{\nu}=\{(x,z)\in\tilde{U}\times B_m(r):
q_i(H_{\nu}(x),z)=0, i=1,\ldots,s\},$$ satis\-fies:\vspace*{2mm}\\
(1) $X_{\nu}$ is $n$-dimensional with proper projection
onto $\tilde{U}$ for almost all $\nu,$\\
(2) $\max\{\sharp(X\cap (\{x\}\times\mathbf{C}^m)):x\in
U\}=$\\
\hspace*{6.0mm}$\max\{\sharp((X_{\nu})_{(n)}\cap
(\{x\}\times\mathbf{C}^m)):x\in \tilde{U}\}$ for almost all $\nu,$\\
(3) $\{X_{\nu}\},$ $\{(X_{\nu})_{(n)}\}$ converge to
$X\cap(\tilde{U}\times\mathbf{C}^m)$ locally
uniformly.\vspace*{1mm}
\end{proposition}

\hspace*{-4.9mm}\textit{Proof of Proposition \ref{ndimensional}.}
Define the algebraic set
$$V=\{(v,z)\in\mathbf{C}^p\times\mathbf{C}^m:q_i(v,z)=0,
i=1,\ldots,s \}.$$

Next, by $F$ denote the intersection of all algebraic subsets of
$\mathbf{C}^p$ containing the image of $H.$ Clearly, $F$ is
irreducible (because $U$ is connected) hence of pure dimension,
say $\bar{n}.$ Fix an open connected subset
$\tilde{U}\subset\subset U.$ In the following lemma $F$ is endowed
with the topology induced by the standard topology of
$\mathbf{C}^{p}.$
%
%*********************************************************
%
\begin{lemma}\label{unbound}Let $r>0$ be such
that $(\tilde{U}\times B_m(r))\cap X\neq\emptyset$ and
$(\overline{\tilde{U}}\times\partial B_m(r))\cap X=\emptyset.$
Then there is an open neighborhood $C$ of
$\overline{H(\tilde{U})}$ in $F$ such that $(C\times B_m(r))\cap
V$ is $\bar{n}$-dimensional with proper projection onto $C.$
Moreover, for every $(a,z)\in(\overline{H(\tilde{U})}\times
B_m(r))\cap V$ it holds $dim_{(a,z)}((C\times B_m(r))\cap
V)=\bar{n}.$
\end{lemma}\vspace*{2mm}

\hspace*{-4.9mm}\textit{Proof of Lemma \ref{unbound}.} First we
check that there is an open neighborhood $C$ of
$\overline{H(\tilde{U})}$ in $F$ such that
$(\overline{C}\times\partial B_m(r))\cap V=\emptyset,$ which
implies the properness of the projection of $(C\times B_m(r))\cap
V$ onto $C.$

It is sufficient to show that for every
$a\in\overline{H(\tilde{U})}$ there is an open neighborhood $C_a$
in $F$ such that $(C_a\times\partial B_m(r))\cap V=\emptyset.$ Fix
$a\in\overline{H(\tilde{U})}.$ Now, if for every open neighborhood
$C_a$ of $a$ we had $(C_a\times\partial B_m(r))\cap
V\neq\emptyset$ then there would be $(\{a\}\times\partial
B_m(r))\cap V\neq\emptyset.$ But then
$(\overline{\tilde{U}}\times\partial B_m(r))\cap X\neq\emptyset$
as $a\in\overline{H(\tilde{U})}\subset H(\overline{\tilde{U}}),$ a
contradiction.

Let us show that $dim_{(a,z)}((C\times B_m(r))\cap V)=\bar{n}$ for
every $(a,z)\in (\overline{H(\tilde{U})}\times B_m(r))\cap V.$
First observe that $dim((C\times B_m(r))\cap V)$ cannot exceed the
dimension of $C$ because $(C\times B_m(r))\cap V$ is with proper
projection onto $C.$ Next suppose that there is $(a,z)\in
(\overline{H(\tilde{U})}\times B_m(r))\cap V$ such that
$dim_{(a,z)}((C\times B_m(r))\cap V)<\bar{n}.$ Let $V_1$ be the
union of the irreducible analytic components of $(C\times
B_m(r))\cap V$ containing $(a,z)$ and let
$\pi:\mathbf{C}^{p}\times\mathbf{C}^m\rightarrow\mathbf{C}^{p}$
denote the natural projection. It is easy to see that
$H^{-1}(\pi(V_1))$ is a non-empty nowhere dense analytic subset of
$H^{-1}(C)$ (nowhere-density because otherwise $H(U)$ would be
contained in an algebraic set of dimension smaller than
$\bar{n}$). Let $P$ be a neighborhood of $(a,z)$ in $C\times
B_m(r)$ such that $P\cap V=P\cap V_1\neq\emptyset.$ Now consider
the set
$$E=\{(w,y)\in (U\times B_m(r))\cap X:(H(w),y)\in P\cap V\}.$$ One
observes that $E\neq\emptyset,$ because
$H^{-1}(\{a\})\times\{z\}\subset E,$ and that $E$ has a non-empty
interior in $X,$ and moreover, the projection of $E$ onto $U$ is
contained in $H^{-1}(\pi(V_1)).$ This contradicts the fact that
$X$ is purely $n$-dimensional.

Since $(\tilde{U}\times B_m(r))\cap X\neq\emptyset$ then
$(H(\tilde{U})\times B_m(r))\cap V\neq\emptyset$ so by what we
have proved so far $(C\times B_m(r))\cap
V$ is $\bar{n}$-dimensional.\qed\\
%
%***********************************************************
%

\hspace*{-5mm}\textit{Proof of Proposition \ref{ndimensional}
(continuation)}. Let $r_0>0$ be such that $(\tilde{U}\times
B_m(r_0))\cap X=(\tilde{U}\times\mathbf{C}^m)\cap X$ and let
$r>r_0.$ Then $(\overline{\tilde{U}}\times\partial B_m(r))\cap
X=\emptyset$ and by Lemma \ref{unbound}, there is a neighborhood
$C$ of $\overline{H(\tilde{U})}$ in $F$ such that $(C\times
B_m(r))\cap V$ is $\bar{n}$-dimensional with proper projection
onto $C.$ Moreover, for every
$(a,z)\in(\overline{H(\tilde{U})}\times B_m(r))\cap V$ it holds
$dim_{(a,z)}((C\times B_m(r))\cap V)=\bar{n}.$ Let
$\{H_{\nu}:\tilde{U}\rightarrow F\}$ be a sequence of holomorphic
mappings converging uniformly to $H$ on $\tilde{U}.$ Define the
sequence $\{X_{\nu}\}$ as in the statement of Proposition
~\ref{ndimensional}.

First we show (1): $X_{\nu}$ is $n$-dimensional and with proper
projection onto $\tilde{U}$ for almost all $\nu.$ To do this
observe that for sufficiently large $\nu$ it holds
$H_{\nu}(\tilde{U})\subset C$ and then
$$X_{\nu}=\{(x,z)\in\tilde{U}\times B_m(r): (H_{\nu}(x),z)\in
(C\times B_m(r))\cap V\}.$$ Thus the properness of the projection
of $X_{\nu}$ onto $\tilde{U}$ is obvious by the choice of $C$ in
Lemma \ref{unbound}.

Now we check the following claim: for sufficiently large $\nu$
every fiber in $X_{\nu}$ over $\tilde{U}$ is not empty. Indeed,
let $C_0$ denote the irreducible Nash component of $C$ containing
$H(\tilde{U}).$ Then the projection of $(C_0\times B_m(r))\cap V$
onto $C_0$ is surjective which follows by Lemma \ref{unbound}. On
the other hand, for sufficiently large $\nu,$
$H_{\nu}(\tilde{U})\subset C_0$ which clearly implies the claim.
Consequently, $X_{\nu}$ is $n$-dimensional for almost all $\nu.$

Let us turn to (2).  Since $C_0$ is an irreducible Nash set then
$Reg(C_0)$ is connected. There is a nowhere dense Nash subset $C'$
of $C_0$ such that the function $\rho: Reg(C_0)\setminus
C'\rightarrow\mathbf{N}$ given by
$$\rho(v)=\sharp((\{v\}\times B_m(r))\cap V)$$
is constant. By $\tilde{m}$ we denote the only value of $\rho.$

Neither $H(\tilde{U})$ nor $H_{\nu}(\tilde{U})$ (for large $\nu$)
can be contained in $Sing(C_0)\cup C'$ so $(H^{-1}(Sing(C_0)\cup
C')\cup H_{\nu}^{-1}(Sing(C_0)\cup C'))\cap \tilde{U}$ is a
nowhere dense analytic subset of $\tilde{U}.$ This means that for
the generic $x\in\tilde{U}$ the fibers in $X$ and in $X_{\nu}$
over $x$ have $\tilde{m}$ elements which completes the proof of
(2).

Finally, let us prove (3). To check the condition (2l) of the
definition of local uniform convergence it is sufficient to show
that for every $(x_0,z_0)\in
(\tilde{U}\times\mathbf{C}^m)\setminus X$ there is a neighborhood
$D$ of $(x_0,z_0)$ in $\tilde{U}\times\mathbf{C}^m$ such that
$D\cap X_{\nu}=\emptyset$ for almost all $\nu.$ This is obvious as
there is $i\in\{1,\ldots,s\}$ such that $q_i(H(x_0),z_0)\neq 0.$
Then $q_i(H_{\nu}(x_0),z_0)\neq 0$ for almost all $\nu$ in some
neighborhood of $(x_0,z_0).$

As for the condition (1l), it suffices to show that for a fixed
$x_0\in\tilde{U}\setminus H^{-1}(Sing(C))$ the sequence
$\{(\{x_0\}\times\mathbf{C}^m)\cap (X_{\nu})_{(n)}\}$ converges to
$(\{x_0\}\times\mathbf{C}^m)\cap X$ locally uniformly. Take
$(x_0,z_0)\in
X\cap(\tilde{U}\times\mathbf{C}^m)=X\cap(\tilde{U}\times B_m(r)).$
Then by Lemma \ref{unbound} it holds $dim_{(H(x_0),z_0)}(C\times
B_m(r))\cap V=dim(C).$ Consequently, (since $H(x_0)\in Reg(C)$ and
$(C\times B_m(r))\cap V$ is with proper projection onto $C$) there
is a sequence $\{z_{\nu}\}$ converging to $z_0$ such that
$dim_{(H_{\nu}(x_0),z_{\nu})}(C\times B_m(r))\cap V=dim(C)$ for
almost all $\nu.$ This implies that for sufficiently large $\nu,$
the image of the projection of every open neighborhood of
$(x_0,z_{\nu})$ in $X_{\nu}$ onto $\tilde{U}$ contains a
neighborhood of $x_0$ in $\tilde{U}.$ Thus
$(x_0,z_{\nu})\in(X_{\nu})_{(n)}$ for almost
all $\nu$ and the proof is complete.\qed\\

\hspace*{-5.2mm}\textit{Proof of Theorem \ref{main} (end).} Let
$F$ denote the intersection of all algebraic subvarieties of
$\mathbf{C}^{p}$ containing $H(U)$ and let $\{H_{\nu}:U\rightarrow
F\}$ be a sequence of holomorphic mappings converging locally
uniformly to $H.$ Define $X_{\nu}$ as in the statement of Theorem
\ref{main}.

It is sufficient to show that for every relatively compact subset
$\tilde{U}$ of $U$ the sequences
$\{X_{\nu}\cap(\tilde{U}\times\mathbf{C}^m)\}$ and
$\{(X_{\nu})_{(n)}\cap(\tilde{U}\times\mathbf{C}^m)\}$ converge to
$X\cap(\tilde{U}\times\mathbf{C}^m)$ locally uniformly and in the
sense of holomorphic chains respectively. Fix
$\tilde{U}\subset\subset U.$ Then by Proposition
\ref{ndimensional} there is $r_0$ such that for every $r>r_0$ the
following hold. $\{X_{\nu}\cap(\tilde{U}\times B_m(r))\}$ and
$\{(X_{\nu})_{(n)}\cap(\tilde{U}\times B_m(r))\}$ converge to
$X\cap(\tilde{U}\times\mathbf{C}^m)$ locally uniformly. Moreover,
for almost all $\nu,$ $X_{\nu}\cap(\tilde{U}\times B_m(r))$ is
$n$-dimensional with proper projection onto $\tilde{U}$ and
$\max\{\sharp(X\cap (\{x\}\times\mathbf{C}^m)):x\in
\tilde{U}\}=\max\{\sharp((X_{\nu})_{(n)}\cap (\{x\}\times
B_m(r))):x\in \tilde{U}\}.$ Thus by Lemma \ref{eqconv} we have:
$\{(X_{\nu})_{(n)}\cap(\tilde{U}\times B_m(r))\}$ converges to
$X\cap(\tilde{U}\times\mathbf{C}^m)$ in the sense of holomorphic
chains. Since $r$ can be taken arbitrarily large we get our claim.\qed\\
%******************************************************************


\begin{thebibliography}{}
\bibitem{An}
Andr\'e, M.: Cinq expos\'es sur la d\'esingularization,
manuscript, \'Ecole Polytechnique F\'ed\'erale de Lausanne, 1992
\bibitem{Ar2}
Artin, M.: \emph{Algebraic approximation of structures over
complete local rings.} Publ. I.H.E.S. \textbf{36}, 23-58 (1969)
\bibitem{Ar3}
Artin, M.: \emph{Algebraic structure of power series rings.}
Contemp. Math. \textbf{13}, 223-227 (1982)
\bibitem{Ba}
Barlet, D.: \textit{Espace analytique r\'eduit des cycles
analytiques complexes compacts d'un espace analytique complexe de
dimension finie.} Fonctions de plusieurs variables complexes, II,
S\'em. Fran\c{c}ois Norguet, 1974-1975, Lecture Notes in Math.,
\textbf{482}, pp. 1-158, Springer, Berlin 1975
\bibitem{B1}
Bilski, M.: \textit{On approximation of analytic sets.}
Manuscripta Math. \textbf{114}, 45-60 (2004)
\bibitem{B2}
Bilski, M.: \textit{Approximation of analytic sets with proper
projection by Nash sets.} C.R. Acad. Sci. Paris, Ser. I
\textbf{341}, 747-750 (2005)
\bibitem{B3}
Bilski, M.: \textit{Approximation of analytic sets by Nash
tangents of higher order.} Math. Z. \textbf{256}, 705-716 (2007)
\bibitem{B4}
Bilski, M.: \textit{Algebraic approximation of analytic sets and
mappings.} Preprint 2007.
\bibitem{BoK}
Bochnak, J., Kucharz, W.: \emph{Approximation of holomorphic maps
by algebraic morphisms.} Ann. Polon. Math. \textbf{80}, 85-92
(2003)
\bibitem{BMT1}
Braun, R. W., Meise, R., Taylor, B. A.: \textit{Higher order
tangents to analytic varieties along curves.} Canad. J. Math.
\textbf{55}, 64-90 (2003)
\bibitem{BMT2}
Braun, R. W., Meise, R., Taylor, B. A.: \textit{The geometry of
analytic varieties satisfying the local Phragm\'en-Lindel\"of
condition and a geometric characterization of the partial
differential operators that are surjective on ${\mathcal
A}(\mathbf{R}^4).$} Trans. Amer. Math. Soc. \textbf{356},
1315-1383 (2004)
\bibitem{Ch}
Chirka, E. M.: Complex analytic sets. Kluwer Academic Publ.,
Dordrecht-Boston-London 1989
\bibitem{DLS}
Demailly, J.-P., Lempert, L., Shiffman, B.: \textit{Algebraic
approximation of holomorphic maps from Stein domains to projective
manifolds.} Duke Math. J. \textbf{76}, 333-363 (1994)
\bibitem{Fo}
Forstneri\v c, F.: \textit{Holomorphic flexibility properties of
complex manifolds.} Amer. J. Math. \textbf{128}, 239-270 (2006)
\bibitem{Ku}
Kucharz, W.: \emph{The Runge approximation problem for holomorphic
maps into Grassmannians.} Math. Z. \textbf{218}, 343-348 (1995)
\bibitem{Lel}
Lelong, P.: \emph{Int\'egration sur un ensemble analytique
complexe.} Bull. Soc. Math. France \textbf{85}, 239-262 (1957)
\bibitem{Lem}
Lempert, L.: \textit{Algebraic approximations in analytic
geometry.} Invent. Math. \textbf{121}, 335-354 (1995)
\bibitem{Og}
Ogoma, T.: \emph{General N\'eron desingularization based on the
idea of Popescu.} J. of Algebra \textbf{167}, 57-84 (1994)
\bibitem{Po1}
Popescu, D.: \emph{General N\'eron desingularization.} Nagoya
Math. J. \textbf{100}, 97-126 (1985)
\bibitem{Po2}
Popescu, D.: \emph{General N\'eron desingularization and
approximation.} Nagoya Math. J. \textbf{104}, 85-115 (1986)
\bibitem{Sp}
Spivakovsky, M.: \emph{A new proof of D. Popescu's theorem on
smoothing of ring homomorphisms.} J. Amer. Math. Soc.,
\textbf{12}, 381-444 (1999)
\bibitem{TT1}
Tancredi, A., Tognoli, A.: \textit{Relative approximation theorems
of Stein manifolds by Nash manifolds.} Boll. Un. Mat. It.
\textbf{3-A,} 343-350 (1989)
\bibitem{TT2}
Tancredi, A., Tognoli, A.: \textit{On the extension of Nash
functions.} Math. Ann. \textbf{288}, 595-604 (1990)
\bibitem{TT3}
Tancredi, A., Tognoli, A.: \textit{On the relative Nash
approximation of analytic maps.} Rev. Mat. Complut. \textbf{11},
185-201 (1998)
\bibitem{Tw}
Tworzewski, P.: \textit{Intersections of analytic sets with linear
subspaces.} Ann. Sc. Norm. Super. Pisa \textbf{17,} 227-271 (1990)
\bibitem{Tw2}
Tworzewski, P.: \textit{Intersection theory in complex analytic
geometry.} Ann. Polon. Math., \textbf{62.2} 177-191 (1995)
\bibitem{Wh}
Whitney, H.: Complex Analytic Varieties. Addison-Wesley Publishing
Co., Reading, Mass.-London-Don Mills, Ont., 1972
\end{thebibliography}
\end{document}